\documentclass[a4paper,11pt]{delfim}

\usepackage{graphicx}

\begin{document}


\title{Remarks on the calculus of variations on time scales\footnote{This work
is part of the first author's PhD project.}}

\author{\textbf{Rui A. C. Ferreira$^1$ and Delfim F. M. Torres$^2$}}

\date{$^1$Department of Mathematics\\
University of Aveiro\\
3810-193 Aveiro, Portugal\\
ruiacferreira@yahoo.com\\ [0.3cm]
$^2$Department of Mathematics\\
University of Aveiro\\
3810-193 Aveiro, Portugal\\
delfim@ua.pt}

\maketitle


\begin{abstract}
\noindent \emph{The calculus of variations is a classical subject
which has gain throughout the last three hundred years a level of
rigor and elegance that only time can give. In this note we show
that, contrary to the classical field, available formulations and
results on the recent calculus of variations on time scales are
still at the heuristic level.}

\medskip


\medskip

\end{abstract}


\section{Introduction}

The calculus on time scales was introduced by Stefan Hilger during
his PhD project, carried out under the scientific supervision of
Bernd Aulbach. The new theory unifies continuous and discrete
analysis and, at the same time, can be applied to other closed
subsets of the real numbers.

The calculus of variations on time scales is in its infancy,
still being possible to give reference to all the works on the
subject: \cite{morian,econo,CD:Bohner:2004,dicovots,comRui:TS:Lisboa07,comRui:TS:LP,zeidan}.

The calculus on time scales has a very similar notation with the
differential calculus. When one is not carefully enough, this
leads to ``results'' without any meaning. This means that extra
care must be taken when dealing with problems on time scales.

In this note we show and highlight some weaknesses that may arise
when proving results on a general time scale. Such weaknesses can
be of various kinds and we hope that by the end of this note they
can be better understood.

Unless the contrary, throughout the text the notation conforms to
that used in the references.


\section{Basics on the time scale calculus}
\label{sec2}

A nonempty closed subset of $\mathbb{R}$ is called a \emph{time
scale} and is denoted by $\mathbb{T}$.

The \emph{forward jump operator}
$\sigma:\mathbb{T}\rightarrow\mathbb{T}$ is defined by
$$\sigma(t)=\inf{\{s\in\mathbb{T}:s>t}\},\mbox{ for all $t\in\mathbb{T}$},$$
while the \emph{backward jump operator}
$\rho:\mathbb{T}\rightarrow\mathbb{T}$ is defined by
$$\rho(t)=\sup{\{s\in\mathbb{T}:s<t}\},\mbox{ for all
$t\in\mathbb{T}$},$$ with $\inf\emptyset=\sup\mathbb{T}$
(\textrm{i.e.}, $\sigma(M)=M$ if $\mathbb{T}$ has a maximum $M$)
and $\sup\emptyset=\inf\mathbb{T}$ (\textrm{i.e.}, $\rho(m)=m$ if
$\mathbb{T}$ has a minimum $m$).

A point $t\in\mathbb{T}$ is called \emph{right-dense},
\emph{right-scattered}, \emph{left-dense} and
\emph{left-scattered} if $\sigma(t)=t$, $\sigma(t)>t$, $\rho(t)=t$
and $\rho(t)<t$, respectively.

We define $\mathbb{T}^k=\mathbb{T}\backslash(\rho(b),b]$ and
$\mathbb{T}^{k^2}=\left(\mathbb{T}^k\right)^k$ and whenever we
write $[a,b]$, we mean $[a,b]\cap\mathbb{T}$ for a given time
scale $\mathbb{T}$.

The \emph{graininess function}
$\mu:\mathbb{T}\rightarrow[0,\infty)$ is defined by
$$\mu(t)=\sigma(t)-t,\mbox{ for all $t\in\mathbb{T}$}.$$

We say that a function $f:\mathbb{T}\rightarrow\mathbb{R}$ is
\emph{delta differentiable} at $t\in\mathbb{T}^k$ if there is a
number $f^{\Delta}(t)$ such that for all $\varepsilon>0$ there
exists a neighborhood $U$ of $t$ (\textrm{i.e.},
$U=(t-\delta,t+\delta)\cap\mathbb{T}$ for some $\delta>0$) such
that
$$|f(\sigma(t))-f(s)-f^{\Delta}(t)(\sigma(t)-s)|
\leq\varepsilon|\sigma(t)-s|,\mbox{ for all $s\in U$}.$$ We call
$f^{\Delta}(t)$ the \emph{delta derivative} of $f$ at $t$.

For delta differentiable $f$ and $g$, the next formulas hold:

\begin{align}
f^\sigma(t)&=f(t)+\mu(t)f^\Delta(t)\label{transfor}\\
(fg)^\Delta(t)&=f^\Delta(t)g^\sigma(t)+f(t)g^\Delta(t)\nonumber\\
&=f^\Delta(t)g(t)+f^\sigma(t)g^\Delta(t)\label{transfor2},
\end{align}
where we abbreviate $f\circ\sigma$ by $f^\sigma$.

A function $f:\mathbb{T}\rightarrow\mathbb{R}$ is called
\emph{rd-continuous} if it is continuous at right-dense points and
if the left-sided limit exists at left-dense points. We denote the
set of all rd-continuous functions by C$_{\textrm{rd}}$ or
C$_{\textrm{rd}}[\mathbb{T}]$, and the set of all delta
differentiable functions with rd-continuous derivative by
C$_{\textrm{rd}}^1$ or C$_{\textrm{rd}}^1[\mathbb{T}]$.

It is known that rd-continuous functions possess an
\emph{antiderivative}, \textrm{i.e.}, there exists a function $F$
with $F^\Delta=f$, and in this case an \emph{integral} is defined
by $\int_{a}^{b}f(t)\Delta t=F(b)-F(a)$. It satisfies
\begin{equation}
\label{sigma} \int_t^{\sigma(t)}f(\tau)\Delta\tau=\mu(t)f(t) \, .
\end{equation}

We now present the integration by parts formulas of the delta
integral:

\begin{lemma}
\label{integracao:partes} If $a,b\in\mathbb{T}$ and
$f,g\in$C$_{\textrm{rd}}^1$, then
\begin{enumerate}

 \item$\int_{a}^{b}f(\sigma(t))g^{\Delta}(t)\Delta t
 =\left[(fg)(t)\right]_{t=a}^{t=b}-\int_{a}^{b}f^{\Delta}(t)g(t)\Delta
 t$;

\item $\int_{a}^{b}f(t)g^{\Delta}(t)\Delta t
=\left[(fg)(t)\right]_{t=a}^{t=b}-\int_{a}^{b}f^{\Delta}(t)g(\sigma(t))\Delta
t$.
\end{enumerate}
\end{lemma}


\section{On the calculus of variations on time scales}

We point out some of the issues which are not completely clear in
the results available in the literature. We do not claim the main
results to be wrong, we just comment on things which are not clear
to us. We hope our remarks to be of some usefulness for the
interested reader on the calculus of variations on time scales,
and we look forward to readers comments and insights.


\subsection{Critical reading of \cite{zeidan}}

We start with some comments on \cite{zeidan}. There, the basic
problem of the calculus of variations on time scales with variable
endpoints, in the class of weak local piecewise $C\sp 1\sb {\rm
rd}$ functions, is studied. Among other things, the following
transversality condition is obtained:
\begin{equation}
\label{transvzeidan} (\hat{L}_v(a),-\hat{L}_v(b))=\nabla
K(\hat{y}(a),\hat{y}(b))+\gamma^T M \, .
\end{equation}
As the authors point out, the delta derivative of a function is
not well defined at a left-scattered maximum point of
$\mathbb{T}$. In virtue of this, since $L$ depends on $y^\Delta$,
the variable $t$ in $L$ is defined ``only'' in $[a,\rho(b)]$. The
Euler-Lagrange equation (in integral form) is
\begin{equation}
\label{e-l} \hat{L}_v(t)=\int_a^t\hat{L}_y(\tau)\Delta\tau+c^T,\
t\in[a,\rho(b)].
\end{equation}
Proof of (\ref{transvzeidan}) in \cite{zeidan} uses the fact that
$\hat{L}_v(b)=\int_a^b\hat{L}_y(\tau)\Delta\tau+c^T$, which is not
included in (\ref{e-l}) if $b$ is a left-scattered point.


\subsection{Critical reading of \cite{econo}}

Now we discuss \cite[Lemma~2.1]{econo}. We first enunciate it:

\begin{lemma} (\cite[Lemma~2.1]{econo})
If $f(t)$ is continuous on $[\rho(a),b]$, where $\rho(a)<b$, and
if
$$\int_{\rho(a)}^b f(t)g(t)\nabla t=0$$
for every function $g(t)\in C[\rho(a),b]$ with
$g(\rho(a))=g(b)=0$, then $f(t)=0$ for $t\in[\rho(a),b]$.
\end{lemma}

We claim that the lemma is not proved for $t=\rho(a)$ and $t=b$.
To see this let us consider the time scale
$\mathbb{T}=\{1,2,3,4,5\}$. Then, every function is continuous.
Let $g$ be an arbitrary function such that $g(1)=g(5)=0$. Define
$f$ in $\mathbb{T}$ by
\[ f(t) = \left\{ \begin{array}{ll}
1 & \mbox{if $t =1$ or $t=5$};\\
0 & \mbox{otherwise}.\end{array} \right. \] Then,
$$\int_{\rho(a)}^b f(t)g(t)\nabla t=\sum_{t=2}^{5}f(t)g(t)=0,$$
but $f(\rho(a))$ and $f(b)$ are not zero.


\subsection{Critical reading of \cite{morian}}

In what follows, we will make some observations about the proof of
the fundamental lemma of the calculus of variations and the
derivation of the Euler-Lagrange equation in \cite{morian}.

The fundamental lemma for the one independent variable problem of
the calculus of variations on time scales, as stated in
\cite{morian}, is (here we do not use the notation of
\cite{morian}):

\begin{lemma} (\cite[Lemma~8]{morian})
If $M$ is a continuous function on $[a,b]^{k^2}$, and if
$$\int_a^b M(t)\eta^\sigma(t)\Delta t=0$$
for all $\eta\in C^1$ with $\eta(a)=\eta(b)=0$, then $M(t)=0$ on
$[a,b]^{k^2}$.
\end{lemma}

We first note that in differential calculus, a $C^1$ function is
(by definition) a function that possess a derivative and this
derivative is continuous. Therefore, it is natural that in the
context of time scales, the notation $C^1$ should mean that a
function has a delta derivative and that the delta derivative is
continuous.

The proof of the lemma is made in various steps (cases).
Essentially, it is made by contradiction, i.e, assuming that
$M(t_0)>0$ for some $t_0\in[a,b]^{k^2}$. It is constructed a
function $\eta$ that is positive in the same points of $M$ and
zero at the other points. In case $A$ of the proof in
\cite{morian}, the authors intend to show that the lemma is true
for left-dense points (independently of the points being
right-dense or right-scattered). The point we want to emphasize
now is that the function constructed in the proof does not
necessarily belong to $C^1$. The above mentioned function is
\[ \eta(t) = \left\{ \begin{array}{ll}
(t-\sigma(u_1))^2(\sigma(t_0)-t)^2 & \mbox{if $t\in[\sigma(u_1),\sigma(t_0)]$},\\
0 & \mbox{otherwise},\end{array} \right. \] with $\sigma(u_1)<t_0$
($\sigma(u_1)\in[t_0-\delta,t_0)$, for some $\delta>0$). We now
assume that $t_0$ is right-scattered. Function $\eta$ is
delta-differentiable in $(\sigma(u_1),\sigma(t_0))$ since it is
the product of two delta differentiable functions, and by
(\ref{transfor2})

\begin{align}
\eta^\Delta(t)&=[(t-\sigma(u_1))^2]^\Delta(\sigma(t_0)-t)^2
+(\sigma(t)-\sigma(u_1))^2[(\sigma(t_0)-t)^2]^\Delta\nonumber\\
&=(t+\sigma(t)-2\sigma(u_1))(\sigma(t_0)-t)^2+(\sigma(t)
-\sigma(u_1))^2(t+\sigma(t)-2\sigma(t_0))\label{rui5}.
\end{align}
What is important to note here is that $\sigma$ is not continuous
at left-dense right-scattered points. This is due to the fact that
if $t$ is such a point, then $$\lim_{s\rightarrow
t^-}\sigma(s)=t<\sigma(t).$$ By (\ref{rui5}), it immediately
follows that $\eta^\Delta$ is not continuous at $t_0$. Clearly,
analogous observations to those above can be made to the double
integral problem.

\bigskip

Now we turn our attention to \cite[page~48]{morian}. There, after
applying the integration by parts formula with respect to $x$, it
appears the following term within a formula (see
\cite[page~48]{morian} for a better comprehension of the
notation):

$$\int_c^d L_p(b,y,\hat{z}(\sigma(b),\tau(y)),
\hat{z}^\Gamma(b,\tau(y)),\hat{z}^\Delta(\sigma(b),y))\zeta(b,y)\Delta
y.$$ By definition, $b=\max\mathbb{X}$, so by the same reasoning
as we did before, we may say that $\hat{z}^\Gamma(b,\tau(y))$
doesn't make sense. However, it is possible to bypass this problem
and we shall show how.

Multiple integration on time scales was introduced in
\cite{miots}. There, it is defined the double Riemann integral. We
use this concept and we define the problem of minimizing the
functional
\begin{equation}
\label{P}
J(u)=\int_{a_1}^{b_1}\int_{a_2}^{b_2}L(t_1,t_2,u(\sigma_1(t_1),
\sigma_2(t_2)),u^{\Delta_1}(t_1,\sigma_2(t_2)),u^{\Delta_2}(\sigma_1(t_1),t_2))\Delta
t_2\Delta t_1,
\end{equation}
among all the functions $u$ that have partial derivatives of the
second order with respect to its arguments, and that become a
given continuous function on the boundary of
$R\doteq[a_1,b_1]\times[a_2,b_2]\subset\mathbb{T}_1\times\mathbb{T}_2$.
We assume that $\sigma_1$ and $\sigma_2$ are delta differentiable.
This immediately implies that $u(\sigma_1(t_1),\sigma_2(t_2))$ is
a continuous function for $(t_1,t_2)\in R$. Further, let
$L(t_1,t_2,y_0,y_1,y_2):\mathbb{T}_1^k\times\mathbb{T}_2^k\times\mathbb{R}^3\rightarrow
\mathbb{R}$ have the necessary smoothness properties in order for
the calculations made below to make sense.

\begin{remark}
If $f(t_1,t_2)$ is continuous in a ``rectangle"\ $R$, then we can
interchange the order of integration in
$$\int\int_R f(t_1,t_2)\Delta t_1\Delta t_2.$$
\end{remark}

We only want to show how to eliminate the problem that appears in
\cite{morian}; so we start with
\begin{equation}
\label{ded0}
\int_{a_1}^{b_1}\int_{a_2}^{b_2}\left[L_{y_0}(\cdot\cdot)\eta(\sigma_1(t_1),
\sigma_2(t_2))+L_{y_1}(\cdot\cdot)\eta^{\Delta_1}(t_1,\sigma_2(t_2))
+L_{y_2}(\cdot\cdot)\eta^{\Delta_2}(\sigma_1(t_1),t_2)\right]\Delta
t_2\Delta t_1,
\end{equation}
where
$(\cdot\cdot)=(t_1,t_2,\tilde{u}(\sigma_1(t_1),\sigma_2(t_2)),
\tilde{u}^{\Delta_1}(t_1,\sigma_2(t_2)),
\tilde{u}^{\Delta_2}(\sigma_1(t_1),t_2))$
and $\eta$ is a function that have partial derivatives of the
second order with respect to its arguments and is zero on the
boundary of $R$. We proceed as follows,
\begin{align}
&\int_{a_1}^{b_1}\int_{a_2}^{b_2}\left[L_{y_0}(\cdot\cdot)\eta(\sigma_1(t_1),
\sigma_2(t_2))+L_{y_1}(\cdot\cdot)\eta^{\Delta_1}(t_1,\sigma_2(t_2))
+L_{y_2}(\cdot\cdot)\eta^{\Delta_2}(\sigma_1(t_1),t_2)\right]\Delta
t_2\Delta t_1 \nonumber\\
&=\int_{a_1}^{b_1}\int_{a_2}^{\rho_2(b_2)}\left[L_{y_0}(\cdot\cdot)\eta(\sigma_1(t_1),
\sigma_2(t_2))+L_{y_1}(\cdot\cdot)\eta^{\Delta_1}(t_1,\sigma_2(t_2))
+L_{y_2}(\cdot\cdot)\eta^{\Delta_2}(\sigma_1(t_1),t_2)\right]\Delta
t_2\Delta t_1\nonumber\\
&\quad+\int_{a_1}^{b_1}\int_{\rho_2(b_2)}^{b_2}\left[L_{y_0}(\cdot\cdot)\eta(\sigma_1(t_1),
\sigma_2(t_2))+L_{y_1}(\cdot\cdot)\eta^{\Delta_1}(t_1,\sigma_2(t_2))
+L_{y_2}(\cdot\cdot)\eta^{\Delta_2}(\sigma_1(t_1),t_2)\right]\Delta
t_2\Delta t_1\nonumber\\
&=\int_{a_2}^{\rho_2(b_2)}\int_{a_1}^{\rho_1(b_1)}\left[L_{y_0}(\cdot\cdot)\eta(\sigma_1(t_1),
\sigma_2(t_2))+L_{y_1}(\cdot\cdot)\eta^{\Delta_1}(t_1,\sigma_2(t_2))
+L_{y_2}(\cdot\cdot)\eta^{\Delta_2}(\sigma_1(t_1),t_2)\right]\Delta t_1\Delta t_2\nonumber\\
&\quad+\int_{a_2}^{\rho_2(b_2)}\int_{\rho_1(b_1)}^{b_1}\left[L_{y_0}(\cdot
\cdot)\eta(\sigma_1(t_1),\sigma_2(t_2))+L_{y_1}(\cdot\cdot)\eta^{\Delta_1}(t_1,
\sigma_2(t_2))+L_{y_2}(\cdot\cdot)\eta^{\Delta_2}(\sigma_1(t_1),
t_2)\right]\Delta t_1\Delta t_2\nonumber\\
&\quad+\int_{a_1}^{b_1}\int_{\rho_2(b_2)}^{b_2}\left[L_{y_0}(\cdot
\cdot)\eta(\sigma_1(t_1),\sigma_2(t_2))+L_{y_1}(\cdot\cdot)\eta^{\Delta_1}(t_1,
\sigma_2(t_2))+L_{y_2}(\cdot\cdot)\eta^{\Delta_2}(\sigma_1(t_1),t_2)\right]\Delta
t_2\Delta t_1.\label{ded1}
\end{align}
The first double integral in the last equality becomes
\begin{align}
&\int_{a_2}^{\rho_2(b_2)}\int_{a_1}^{\rho_1(b_1)}\left[L_{y_0}(\cdot
\cdot)\eta(\sigma_1(t_1),\sigma_2(t_2))+L_{y_1}(\cdot
\cdot)\eta^{\Delta_1}(t_1,\sigma_2(t_2))+L_{y_2}(\cdot
\cdot)\eta^{\Delta_2}(\sigma_1(t_1),t_2)\right]\Delta
t_1\Delta t_2\nonumber\\
&=\int_{a_2}^{\rho_2(b_2)}\int_{a_1}^{\rho_1(b_1)}\left[L_{y_0}(\cdot
\cdot)\eta(\sigma_1(t_1),\sigma_2(t_2))-L_{y_1}^{\Delta_1}(\cdot
\cdot)\eta(\sigma_1(t_1),\sigma_2(t_2))-L_{y_2}^{\Delta_2}(\cdot
\cdot)\eta(\sigma_1(t_1),\sigma_2(t_2))\right]
\Delta t_1\Delta t_2\nonumber\\
&\quad+\int_{a_2}^{\rho_2(b_2)}\int_{a_1}^{\rho_1(b_1)}\left\{
[L_{y_1}(\cdot\cdot)\eta(t_1,\sigma_2(t_2))]^{\Delta_1}
+[L_{y_2}(\cdot\cdot)\eta(\sigma_1(t_1),t_2]^{\Delta_2}\right\}\Delta
t_1\Delta t_2\nonumber\\
&=\int_{a_2}^{\rho_2(b_2)}\int_{a_1}^{\rho_1(b_1)}\left[L_{y_0}(\cdot
\cdot)-L_{y_1}^{\Delta_1}(\cdot\cdot)-L_{y_2}^{\Delta_2}(\cdot
\cdot)\right]\eta(\sigma_1(t_1),\sigma_2(t_2))
\Delta t_1\Delta t_2\nonumber \\
&\quad+\int_{a_2}^{\rho_2(b_2)}\left[L_{y_1}(\cdot\cdot)\eta(t_1,
\sigma_2(t_2))\right]_{t_1=a_1}^{t_1=\rho_1(b_1)}\Delta
t_2+\int_{a_1}^{\rho_1(b_1)}\left[L_{y_2}(\cdot\cdot)\eta(\sigma_1(t_1),
t_2)\right]_{t_2=a_2}^{t_2=\rho_2(b_2)}\Delta
t_1\nonumber\\
&=\int_{a_2}^{\rho_2(b_2)}\int_{a_1}^{\rho_1(b_1)}\left[L_{y_0}(\cdot
\cdot)-L_{y_1}^{\Delta_1}(\cdot\cdot)-L_{y_2}^{\Delta_2}(\cdot
\cdot)\right]\eta(\sigma_1(t_1),\sigma_2(t_2))
\Delta t_1\Delta t_2\nonumber \\
&\quad+\int_{a_2}^{\rho_2(b_2)}L_{y_1}(\rho_1(b_1)\cdot)
\eta(\rho_1(b_1),\sigma_2(t_2))\Delta
t_2+\int_{a_1}^{\rho_1(b_1)}L_{y_2}(\cdot
\rho_2(b_2))\eta(\sigma_1(t_1),\rho_2(b_2))\Delta
t_1\label{ded5},
\end{align}
where $(\rho_1(b_1)\cdot)$, $(\cdot\rho_2(t_2))$ is equal to
$(\cdot\cdot)$ with $t_1=\rho_1(b_1)$, $t_2=\rho_2(b_2)$,
respectively.

The second double integral in (\ref{ded1}) becomes
\begin{align}
&\int_{a_2}^{\rho_2(b_2)}\int_{\rho_1(b_1)}^{b_1}\left[L_{y_0}(\cdot
\cdot)\eta(\sigma_1(t_1),\sigma_2(t_2))+L_{y_1}(\cdot\cdot)\eta^{\Delta_1}(t_1,
\sigma_2(t_2))+L_{y_2}(\cdot\cdot)\eta^{\Delta_2}(\sigma_1(t_1),t_2)\right]\Delta
t_1\Delta t_2\nonumber\\
&=\int_{a_2}^{\rho_2(b_2)}\mu_1(\rho_1(b_1))\left\{L_{y_1}(\rho_1(b_1)\cdot)
\eta^{\Delta_1}(\rho_1(b_1),\sigma_2(t_2))+L_{y_2}(\rho_1(b_1)
\cdot)\eta^{\Delta_2}(\sigma_1(\rho_1(b_1)),t_2)\right\}\Delta
t_2\label{ded2},
\end{align}
because
$L_{y_0}(\rho_1(b_1)\cdot)\eta(\sigma_1(\rho_1(b_1)),
\sigma_2(t_2))=L_{y_0}(\rho_1(b_1)\cdot)\eta(b_1,\sigma_2(t_2))=0$.
Now, note that
$$\eta(b_1,\sigma_2(t_2))=\eta(\rho_1(b_1),\sigma_2(t_2))
+\mu_1(\rho_1(b_1))\eta^{\Delta_1}(\rho_1(b_1),\sigma_2(t_2)),$$
hence
$\mu_1(\rho_1(b_1))\eta^{\Delta_1}(\rho_1(b_1),\sigma_2(t_2))
=-\eta(\rho_1(b_1),\sigma_2(t_2))$. Therefore, (\ref{ded2})
becomes
\begin{equation}
\label{ded3}
\int_{a_2}^{\rho_2(b_2)}\left\{-L_{y_1}(\rho_1(b_1)\cdot)
\eta(\rho_1(b_1),\sigma_2(t_2))+\mu_1(\rho_1(b_1))L_{y_2}(\rho_1(b_1)\cdot)
\eta^{\Delta_2}(\sigma_1(\rho_1(b_1)),t_2)\right\}\Delta t_2.
\end{equation}
Now,
\begin{align}
&\int_{a_2}^{\rho_2(b_2)}L_{y_2}(\rho_1(b_1)\cdot)
\eta^{\Delta_2}(\sigma_1(\rho_1(b_1)),t_2)\Delta t_2\nonumber\\
&=
\int_{a_2}^{\rho_2(b_2)}L_{y_2}(\rho_1(b_1)\cdot)\eta^{\Delta_2}(b_1,t_2)\Delta
t_2\nonumber\\
&=\left[L_{y_2}(\rho_1(b_1)\cdot)\eta(b_1,t_2)\right]_{t_2=a_2}^{t_2=\rho_2(b_2)}
-\int_{a_2}^{\rho_2(b_2)}L_{y_2}^{\Delta_2}(\rho_1(b_1)\cdot)\eta(b_1,
\sigma_2(t_2))\Delta t_2\nonumber\\
&=0\nonumber,
\end{align}
hence (\ref{ded3}) becomes
\begin{equation}
\label{ded4}
\int_{a_2}^{\rho_2(b_2)}-L_{y_1}(\rho_1(b_1)\cdot)\eta(\rho_1(b_1),\sigma_2(t_2))\Delta
t_2.
\end{equation}
Finally, the third integral in (\ref{ded1}) becomes (repeating
analogous steps as above)
\begin{align}
&\int_{a_1}^{b_1}\int_{\rho_2(b_2)}^{b_2}\left[L_{y_0}(\cdot\cdot)\eta(\sigma_1(t_1),
\sigma_2(t_2))+L_{y_1}(\cdot\cdot)\eta^{\Delta_1}(t_1,\sigma_2(t_2))
+L_{y_2}(\cdot\cdot)\eta^{\Delta_2}(\sigma_1(t_1),t_2)\right]\Delta t_2\Delta t_1\nonumber\\
&=\int_{a_1}^{b_1}\mu_2(\rho_2(b_2))\left[L_{y_1}(\cdot\rho_2(b_2))\eta^{\Delta_1}(t_1,b_2)
+L_{y_2}(\cdot\rho_2(b_2))\eta^{\Delta_2}(\sigma_1(t_1),\rho_2(b_2))\right]\Delta t_1\nonumber\\
&=\int_{a_1}^{\rho_1(b_1)}\mu_2(\rho_2(b_2))\left[L_{y_1}(\cdot\rho_2(b_2))
\eta^{\Delta_1}(t_1,b_2)+L_{y_2}(\cdot\rho_2(b_2))
\eta^{\Delta_2}(\sigma_1(t_1),\rho_2(b_2))\right]\Delta t_1\nonumber\\
&\quad+\mu_1(\rho_1(b_1))\mu_2(\rho_2(b_2))\left[L_{y_1}(\rho_1(b_1),
\rho_2(b_2))\eta^{\Delta_1}(\rho_1(b_1),b_2)+L_{y_2}(\rho_1(b_1),
\rho_2(b_2))\eta^{\Delta_2}(b_1,\rho_2(b_2))\right]\nonumber\\
&=\int_{a_1}^{\rho_1(b_1)}-L_{y_2}(\cdot\rho_2(t_2))\eta(\sigma_1(t_1),
\rho_2(t_2))\Delta t_1\label{ded6},
\end{align}
where $(\rho_1(b_1),\rho_2(b_2))$ are the arguments of $(t_1,t_2)$
on $(\cdot\cdot)$. Combining (\ref{ded0}), (\ref{ded1}),
(\ref{ded5}), (\ref{ded4}) and (\ref{ded6}), we obtain
\begin{equation}
\label{dedf}
\int_{a_2}^{\rho_2(b_2)}\int_{a_1}^{\rho_1(b_1)}\left[L_{y_0}(\cdot\cdot)
-L_{y_1}^{\Delta_1}(\cdot\cdot)-L_{y_2}^{\Delta_2}(\cdot\cdot)\right]
\eta(\sigma_1(t_1),\sigma_2(t_2))
\Delta t_1\Delta t_2.
\end{equation}

\begin{remark}
A version of the fundamental lemma of the calculus of variations
on time scales can be proved in order to deduce the Euler-Lagrange
equation from (\ref{dedf}). Note that, since we assumed that
$\sigma_i\ (i=1,2)$ are $\Delta_i$-differentiable, it is not
possible to have left-dense, right-scattered points at the same
time within the time scales $\mathbb{T}_i$.
\end{remark}


\subsection{Critical reading of \cite{dicovots}}

Finally, we want to make some observations regarding the paper
\cite{dicovots}.

\begin{lemma}(\cite[Lemma~5.1]{dicovots})
If $M(x,y)$ is continuous on $E\cup \Gamma$ with
$$\int\int_E M(x,y)\eta (\sigma_1(x),\sigma_2(y))\Delta_1x\Delta_2y=0$$
for every admissible variation $\eta$, then
$$M(x,y)=0,\ \mbox{for all}\ (x,y)\in E^0.$$
\end{lemma}
We claim that this lemma is not proved. To see this, consider
$$E=[0,5]\cap\mathbb{Z}\times[0,5]\cap\mathbb{Z}.$$ Suppose that we want to
prove the lemma for the point $(1,1)\in E^0$. Continuity of $M$
ensures that $M$ is positive in
$\Omega=[1,\sigma_1(1))\times[1,\sigma_2(1))$. To get a
contradiction, the authors created the function
\[ \eta(x,y) = \left\{ \begin{array}{ll}
(x-x_0)^2[x-\sigma_1(x_1)]^2(y-y_0)^2[y-\sigma_2(y_1)]^2 & \mbox{for $(x,y)\in\Omega$},\\
0 & \mbox{for $(x,y)\in E\backslash\Omega$},\end{array} \right. \]
and concluded that
$$\int\int_E M(x,y)\eta (\sigma_1(x),\sigma_2(y))\Delta_1x\Delta_2y
=\int\int_\Omega M(x,y)\eta
(\sigma_1(x),\sigma_2(y))\Delta_1x\Delta_2y>0.$$ With the notation
used, we have $x_0=y_0=1$. Now note that the set $\Omega$ consists
of only the point $(x_0,y_0)$, and hence $\eta$ is zero in all its
domain. Therefore,
$$\int\int_E M(x,y)\eta (\sigma_1(x),\sigma_2(y))\Delta_1x\Delta_2y=0,$$
which is not a contradiction.

\bigskip

There is another point that we would like to mention. The authors
define an \emph{admissible variation} to be a function in
C$_{rd}^{(1)}(E\cup\Gamma)$ that satisfies $\eta=0$ on $\Gamma$.
The set C$_{rd}^{(1)}$ consists of all continuous functions for
which both the $\Delta_1$-partial derivative and the
$\Delta_2$-partial derivative exist and are of class C$_{rd}$.

When the authors are deducing the Euler-Lagrange equation, a
crucial step is the application of Green's formula to the
expression
$$\int\int_E \left\{\frac{\partial}{\Delta_1 x}[L_p(\cdot)\eta(x,\sigma_2(y))]
+\frac{\partial}{\Delta_2 y}[L_q(\cdot)\eta(\sigma_1(x),y)]
\right\}\Delta_1 x\Delta_2y.$$ In order to apply Green's formula,
we must have the partial delta derivatives in the previous
integral, continuous (see \cite[Theorem~2.25]{dicovots}). If
$\eta\in$ C$_{rd}^{(1)}(E\cup\Gamma)$, why should those partial
derivatives be continuous? A justification to this fact could be
that the authors, when enunciating Theorem~5.2 (Euler's necessary
condition), assume that the admissible functions have continuous
partial delta derivatives of the second order, hence it would be
implicit that the admissible variations should belong to the same
class. However, if this would be the case, why can Lemma~5.1 be
applied to formula
\begin{equation}
\label{alemma} \int\int_E
\left\{L_u(\cdot)-\frac{\partial}{\Delta_1
x}L_p(\cdot)-\frac{\partial}{\Delta_2
y}L_q(\cdot)\right\}\eta(\sigma_1(x),\sigma_2(y))\Delta_1
x\Delta_2y=0\ ?
\end{equation}
It is well know that if a function $f$ is continuous then $f\in$
C$_{rd}$, but the converse is not necessarily true. So, in order
to apply (a version of) Lemma~5.1 to (\ref{alemma}) it would be
necessary to prove it for the class of functions $\eta$ that have
continuous partial delta derivatives of the second order. In
particular, Lemma~5.1 would be a corollary from this previous one.

\bigskip

More could be said, in particular about our own work
\cite{comRui:TS:LP} which is being rewritten in order to
make it correct!


\section*{Acknowledgments}

Work supported by {\it Centre for Research on Optimization and Control} (CEOC)
of the University of Aveiro, through the Portuguese Foundation
for Science and Technology (FCT),
cofinanced by the European Community Fund FEDER/POCI 2010.
The authors are grateful to Calvin D. Ahlbrandt for his
availability in responding to e-mails.



\end{document}